\documentclass[12pt]{article}
\usepackage{amscd,amsmath, amssymb, fancyhdr, mathbbol}

\usepackage{color}

\numberwithin{equation}{section}

\newcommand{\version}{version 1.0,\ \ Nov 5, 2025}

\def\eqref#1{(\ref{#1})}

\newcommand{\arrow}{{\:\longrightarrow\:}}
\newcommand{\Z}{{\mathbb Z}}
\def\C{{\mathbb C}}
\def\P{{\mathbb P}}

\newcommand{\Q}{{\mathbb Q}}

\def\1{\sqrt{-1}\:}

\newcommand{\cntrct}                
{\hspace{2pt}\raisebox{1pt}{\text{$\lrcorner$}}\hspace{2pt}}

\renewcommand{\tilde}{\widetilde}

\renewcommand{\phi}{\varphi}
\renewcommand{\epsilon}{\varepsilon}
\renewcommand{\geq}{\geqslant}

\newcommand{\Teich}{\operatorname{Teich}}

\newcommand{\Pic}{\operatorname{Pic}}

\newcommand{\Sym}{\operatorname{Sym}}

\newcommand{\Diff}{\operatorname{Diff}}

\newcommand{\Hilb}{\operatorname{Hilb}}

\newcommand{\Comp}{\operatorname{Comp}}

\newcounter{Mycounter}[section]
\newcounter{lemma}[section]
\newcounter{claim}[section]
\newcounter{sublemma}[section]
\newcounter{corollary}[section]
\newcounter{theorem}[section]
\renewcommand{\thetheorem}{{Theorem \thesection.\arabic{theorem}}}
\newcommand{\theorem}{%
    \setcounter{theorem}{\value{Mycounter}}
    \refstepcounter{theorem}
    \stepcounter{Mycounter}
    {\noindent \bf \thetheorem:\ }}

\newcounter{conjecture}[section]
\renewcommand{\theconjecture}{{Conjecture \thesection.\arabic{conjecture}}}
\newcommand{\conjecture}{%
    \setcounter{conjecture}{\value{Mycounter}}
    \refstepcounter{conjecture}
    \stepcounter{Mycounter}
    {\noindent \bf \theconjecture:\ }}

\newcounter{proposition}[section]
\renewcommand{\theproposition}
      {{Proposition \thesection.\arabic{proposition}}}
\newcommand{\proposition}{%
    \setcounter{proposition}{\value{Mycounter}}
    \refstepcounter{proposition}
    \stepcounter{Mycounter}
    {\noindent \bf \theproposition:\ }}

\newcounter{definition}[section]
\renewcommand{\thedefinition}
      {{Definition~\thesection.\arabic{definition}}}
\newcommand{\definition}{%
    \setcounter{definition}{\value{Mycounter}}
    \refstepcounter{definition}
    \stepcounter{Mycounter}
    {\noindent \bf \thedefinition:\ }}

\newcounter{example}[section]
\renewcommand{\theexample}{{Example \thesection.\arabic{example}}}
\newcommand{\example}{%
    \setcounter{example}{\value{Mycounter}}
    \refstepcounter{example}
    \stepcounter{Mycounter}
    {\noindent \bf \theexample:\ }}

\newcounter{remark}[section]
\renewcommand{\theremark}{{Remark \thesection.\arabic{remark}}}
\newcommand{\remark}{%
    \setcounter{remark}{\value{Mycounter}}
    \refstepcounter{remark}
    \stepcounter{Mycounter}
    {\noindent \bf \theremark:\ }}

\newcounter{problem}[section]
\newcounter{question}[section]
\newcommand{\proof}{\noindent{\bf Proof:\ }}

\makeatletter

\@addtoreset{equation}{section} \@addtoreset{footnote}{section}
\makeatother

\def\blacksquare{\hbox{\vrule width 5pt height 5pt depth 0pt}}
\def\endproof{\blacksquare}

\begin{document}
\begin{center}
{\LARGE\bf
On anti-hyperbolicity for hyperk\"ahler varieties\\[4mm]
}

Ljudmila Kamenova\footnote{Partially supported 
by a grant from the Simons Foundation/SFARI (522730, LK).}, 
Steven Lu\footnote{partially supported by an NSERC grant.}

\end{center}

{\small \hspace{0.10\linewidth}
\begin{minipage}[t]{0.85\linewidth}
{\bf Abstract} \\
By restricting to (a linear subspace of) an affine chart in projective space, a complex stably rational or unirational manifold of dimension $m$ is meromorphically dominable by $\C^m$, i.e., admits a meromorphic dominating map from $\C^m$. So are varieties that are birational to abelian varieties and Kummer K3 surfaces. G. Buzzard and the second author have shown that elliptic K3 surfaces are holomorphically dominable by $\C^m$, i.e. admitting a holomorphic map with nontrivial jacobian. In this paper we explore various examples and criteria for meromorphic and holomorphic dominability by $\C^m$ of certain hyperk\"ahler manifolds, generalizing some known results about K3 surfaces. 
Anti-hyperbolicity has several interpretations in the sense of vanishing of the Kobayashi-Royden metrics, admitting dense entire holomorphic curves,  or dominating holomorphic or meromorphic maps from the complex affine space of the same dimension. 

\end{minipage}
}

\tableofcontents


\section{Introduction}

It has been conjectured by Green-Griffiths as a natural generalization of 
the Little Picard theorem that a nontrivial holomorphic image of $\C$ in 
a (complex) variety of general type is algebraically degenerate, i.e. not 
Zariski dense. If true, then so would any holomorphic curve be in a variety 
that fibers over a variety of general type. In view of the existence of 
the core map, see \cite{Camp04, Lu02}, any projective variety $X$ admits a meromorphic fibration over a nontrivial ``orbivariety'' base of general type unless $X$ is special. Clearly, an ``orbifold'' generalization of the above conjecture would imply the same algebraic degeneracy condition for holomorphic curves in non-special $X$. So it is interesting to investigate the possible existence of Zariski dense holomorphic curves in special varieties in each dimension $n$, and thus the meromorphic dominability by $\C^n$ of such varieties.

\hfill

Although dominability implies the existence of Zariski dense holomorphic 
curves, the converse is completely open. Easily established dominable 
manifolds include  meromorphic images of abelian 
varieties (such as Kummer varieties), stably rational and unirational manifolds. More recently elliptic K3 surfaces 
have been shown to be dominable \cite{_Bu_Lu_}. These are special varieties in the above sense and all manifolds with trivial canonical bundle such as abelian varieties and K3 surfaces are special. The purpose of the present paper is to generalize the results of \cite{_Bu_Lu_}, which restricts to compact complex surfaces, to a higher dimensional generalization of K3 
surfaces called hyperk\"ahler manifolds or holomorphic symplectic varieties. 
We prove that Hilbert schemes of points on dominable K3 surfaces are 
dominable. The generalized Kummer varieties are also dominable. So are 
hyperk\"ahler manifolds that admit more than one fibration structures.

\hfill

\theorem
Let $M$ be a hyperk\"ahler manifold or a primitive symplectic variety of dimension $2n$. If $M$ admits two distinct Lagrangian fibration structures (necessarily generically transversal), then $M$ is meromorphically dominable by $\C^{2n}$. 

\hfill

In order to prove dominability of hyperk\"ahler manifolds with one 
fibration structure, we will assume in this paper the hoped for fact that the base of the fibration is projective space and that there are no multiple fibers in codimension one of the base. Locally, multiple fibers might exist, 
but globally it is not known if there are examples of multiple fibers 
in codimension one of hyperk\"ahler Lagrangian fibrations. It is a 
classical result that elliptic K3 surfaces do not have multiple fibers 
(see for example, Proposition 1.6. in Ch.11 of \cite{H_K3}). 

\hfill

\theorem
Let $M$ be a hyperk\"ahler manifold or a primitive symplectic variety of dimension $2n$. If $M$ is projective and  admits a Lagrangian fibration $f: M \arrow \C\P^n$ 
with no multiple fibers in codimension one of the base, then $M$ is holomorphically dominable by $\C^{2n}$. 

\hfill 









One of the biggest recent breakthroughs in anti-hyperbolicity since the case of surfaces done almost a quarter century ago in \cite{_Bu_Lu_, BL1, LuMM} is found in the paper of Campana-Winkelmann \cite{CW} where they proved that rationally connected varieties contain dense entire curves. Here we provide natural generalizations of some their other results in special cases to dominability. Given an abelian fibration with no multiple fibers in codimension one over a log-Fano variety, we prove that it contains a dense entire curve. This generalizes Campana-Winkelmann's result \cite[Proposition 7.8]{CW} for elliptic fibrations. 

\hfill

\theorem 
Given a projective manifold $M$ with an abelian fibration $f: M \arrow B$ with no multiple fibres in codimension one over a log-Fano variety $B$, then there exists a dense entire curve.



\section{Basic hyperk\"ahler geometry}


\definition  A {\em hyperk\"ahler manifold} $M$ is a 
compact K\"ahler manifold with $\pi_1(M)=0$ and 
$H^0(M, \Omega^2_M) = \C \omega$, where $\omega$ is an 
everywhere non-degenerate $2$-form.

\hfill

\example
We recall Beauville's standard series of examples, \cite{Beauville}. 
In each possible complex dimension $2n$ there are at least two examples, 
together with their deformations. One of them is the Hilbert scheme 
$\Hilb^n (S)$ of $n$ points on a K3 surface $S$. It can be realized as the 
desingularization of the symmetric product $\Sym^n(S)$ given by the 
Hilbert-Chow morphism $\Hilb^n(S) \arrow \Sym^n(S)$. 
The second example is the 
generalized Kummer variety, $K^n(A)$, of dimension $2n$, 
where $A$ is an abelian surface. Here is 
the basic construction. Consider the Hilbert-Chow morphism 
composed with the summation map given by the group structure on $A$: 
$$\Hilb^{n+1} (A) \arrow \Sym^{n+1} (A) \arrow A.$$ 
Beauville (\cite{Beauville}) showed that the fiber $K^n(A)$ of the 
composite map over $0 \in A$ is a hyperk\"ahler manifold of 
dimension $2n$. 

\hfill

\example
There are two exceptional examples ($O_6$ and $O_{10}$) due to O'Grady 
(\cite{og1,og2}) in complex dimensions $6$ and $10$. 
Let $S$ be a K3 surface, then $O_{10}$ 
is the minimal desingularization of the Mukai moduli space of rank-two 
semistable torsion-free sheaves on $S$ with $c_1=0, c_2=4$. Let $A$ be 
an abelian surface and $\cal M$ be the desingularization of the Mukai moduli 
space of rank-two semistable torsion-free sheaves on $A$ with $c_1=0, c_2=2$. 
Then $O_6$ is a fiber of the natural map 
${\cal M} \rightarrow A \times A^\vee$. 

\hfill

All of the known hyperk\"ahler examples admit fibration structures in a small deformation. By Matsushita's theorem, the possible fibration structure on a hyperk\"ahler manifold is quite restricted. 

\hfill

\theorem (Matsushita, \cite{Matsushita}) 
Let $M$ be a hyperk\"ahler manifold of dimension $2n$, let $B$ be a normal variety with $0 < \dim B < 2n$, and $f: X \arrow B$ be a fibration, i.e., a proper morphism with connected fibers. Then the base $B$ is $n$-dimensional Fano variety with at worst $\Q$-factorial log-terminal singularities. The general fiber of $f$ is an abelian variety which is Lagrangian with respect to the holomorphic symplectic form on $M$. 

\hfill

\remark \label{smooth_base}
In all of the known cases of Lagrangian fibrations on hyperk\"ahler manifolds, the base $B$ is smooth. Hwang \cite{Hwang} has shown that if $M$ is projective and $B$ is smooth, then $B$ is isomorphic to the complex projective space $\C \P^n$. Greb and C. Lehn \cite{GL} generalized this result to the case of non-projective $M$. 
It is conjectured that $B$ is always smooth, and therefore $B \cong \C \P^n$. This conjecture was proven in the four-dimensional case ($n=2$). Ou \cite{Ou} showed that $B$ is isomorphic to either $\C \P^2$ or to a Fano surface with exactly one singular point which is a Du Val singularity of type $E_8$. Huybrechts and Xu \cite{HX} excluded the singular case, and therefore $B \cong \C \P^2$ when $n=2$.

\hfill

Let $M$ be a hyperk\"ahler manifold. 
On $H^2(M, \Z)$ there is a natural primitive integral quadratic form, 
called the {\em Beauville-Bogomolov-Fujiki form}, or the {\em BBF form} 
for short. The simplest way to define it is via the Fujiki relation 
below. The classical definition can be found in \cite{Beauville} and 
\cite{_Huybrechts:basic_}, for example. 

\hfill

\theorem
(Fujiki, \cite{_Fujiki:HK_}) \label{Fujiki_formula}
Let $\eta\in H^2(M, \Z)$ and $\dim M=2n$, where $M$ is a 
hyperk\"ahler manifold. Then $\int_M \eta^{2n}= c \cdot q(\eta,\eta)^n$,
for a primitive integral quadratic form $q$ on $H^2(M, \Z)$, where $c>0$ is a 
constant depending on the topological type of $M$. The constant $c$ in Fujiki's 
formula is called the {\em Fujiki constant}. 
More generally, if $\eta_1, \dots, \eta_{2n} \in H^2(M)$, then 
$$\int_M \eta_1 \wedge \dots \wedge \eta_{2n} = \frac{c}{(2n)!} 
\sum_{\sigma} q(\eta_{\sigma_1}, \eta_{\sigma_2}) \dots 
q(\eta_{\sigma_{2n-1}}, \eta_{\sigma_{2n}}).$$
\endproof

\hfill



\remark \label{proj}
If there is a multisection of a Lagrangian fibration, then the total space is projective. For our purpose, this says that if a hyperk\"ahler manifold admits two Lagrangian fibrations, then the manifold is projective, which we can easily see by taking the sum of two nef classes $L_1, L_2$ defining the two fibrations, and noticing that $q(L_1+L_2) = 2q(L_1, L_2)>0$ (as in the proof of \cite[Theorem 2.11]{klv}), then apply Huybrechts' projectivity criterion. 

\hfill

\definition
Let $\eta \in H^{1,1}(M)$ be a real (1,1)-class on 
a hyperk\"ahler manifold $M$. We say that $\eta$  
is {\em parabolic} if $q(\eta , \eta)=0$. 
A line bundle $L$ is called {\em parabolic} if $c_1(L)$
is parabolic.

\hfill

If $f: M \arrow B$ is a Lagrangian fibration and the class $\alpha \in H^2(B, \Z)$ is ample, then the pull-back $\eta = f^* \alpha$ is nef and trivially satisfies $q(\eta , \eta)=0$ by \ref{Fujiki_formula}. In other words, given a Lagrangian fibration, there is a natural nef parabolic class associated to the fibration. 

\hfill

\definition
Let $L$ be a holomorphic line bundle
on a hyperk\"ahler manifold. We call
$L$ {\em Lagrangian} if it is parabolic and semiample.

\hfill

\conjecture (Hyperk\"ahler SYZ Conjecture) 
If a hyperk\"ahler manifold $M$ has a non-trivial parabolic nef line bundle $L$, then $|mL|$ is base point free for some $m\in \mathbb N$ and defines a Lagrangian fibration $M\to B$. 

\hfill

\definition
Let $M$ be a hyperk\"ahler manifold.
Fix a parabolic class $L\in H^2(M,\Z)$.
We denote by $\Teich_L$ the Teichm\"uller space
of all complex structures $I$ of hyperk\"ahler type 
on $M$ such that $L$ is of type $(1,1)$ on $(M,I)$.
The space $\Teich_L$ is called
{\em the Teichm\"uller space of 
hyperk\"ahler manifolds with parabolic class} and it is a divisor in the 
whole Teichm\"uller space of $M$.

\hfill

\theorem (Kamenova, Verbitsky, \cite{_Kamenova_Verbitsky:fibrations_}) 
\label{_dense_ope_Lagra_Theorem_} 
Let $\Teich_L^\circ\subset\Teich_L$ be the set of all $I\in \Teich_L$ for which
$L$ is Lagrangian. Then  $\Teich_L^\circ$ is dense and open in $\Teich_L$.

\hfill

\definition Let $M$ be a compact complex manifold, and let $\Diff^0(M)$ be a connected component of its diffeomorphism 
group. Denote by $\Comp$ the space of complex structures on $M$, 
equipped with the structure of a Fr\'echet manifold. The 
{\em Teichm\"uller space} of $M$ 
is the quotient  $\Teich:=\Comp/\Diff^0(M)$. 
The Teichm\"uller space is finite-dimensional for a Calabi-Yau $M$, 
see \cite{_Catanese:moduli_}. 

\hfill

\definition
Let $\Diff^+(M)$ be the group of orientable diffeomorphisms of 
a complex manifold $M$. The {\em mapping class group} 
$\Gamma:=\Diff^+(M)/\Diff^0(M)$ acts on $\Teich$. A complex structure 
$I \in \Teich$ is called {\em ergodic} if its orbit $\Gamma \cdot I$ 
is dense in $\Teich$. 

\hfill

\remark \label{rm_double}
From the properties of ergodic structures established by M. Verbitsky 
in \cite{_Verbitsky:ergodic_}, any family given by intrinsic properties is 
dense if it contains an ergodic complex structure. This implies 
that hyperk\"ahler manifolds admitting double fibration structures 
are dense in the Teichm\"uller space (see also \cite{klv}). 

\hfill

A lot of the fundamental results in hyperk\"ahler geometry can be generalized in the singular setting for primitive symplectic varieties. 

\hfill

\definition
A {\em primitive symplectic variety} is a normal compact K\"ahler variety $X$ with rational singularities, such that $H^1(X,\mathcal{O}_X)=0$, and for a symplectic form $\sigma$, 
$H^0(X,\Omega_X^{[2]})=\C\sigma$.

\hfill 

\definition
Let $X$ be a primitive symplectic variety. A {\em rational Lagrangian fibration} is a meromorphic map $f:X \dashrightarrow B$ to a normal K\"ahler variety $B$ such that $f$ has connected fibers and its general fiber is a Lagrangian subvariety of $X$. If $f$ is given by the linear system $|mL|$ for some $L\in Pic(X)$ and $m>>0$  we say that the rational fibration is induced by $L$.

\hfill

\conjecture (Rational SYZ Conjecture) 
Let $X$ be a primitive symplectic variety. If a nontrivial movable line bundle $L$ on $X$ satisfies $q_X(L)=0$, then $L$ induces a rational Lagrangian fibration $f:X \dashrightarrow B$.

\hfill

Let us summarize Matsushita's and Hwang's theorems in the singular setting, as in \cite[Theorem 2.8]{KL}. 

\hfill

\theorem
Let $X$ be a primitive symplectic variety of dimension $2n$ and  let $f: X \rightarrow B$ be a surjective holomorphic map with connected fibers to a normal K\"ahler variety $B$ with $0<\dim B < 2n$. Then the following holds.
\begin{enumerate}
    \item The base $B$ is a projective variety with Picard rank $\rho(B)=1$, in particular, $B$ is projective and has $\Q$-factorial, log-terminal singularities. Furthermore, $\dim B=n$.
    \item The morphism $f$ is equidimensional and each irreducible component of each fiber of $f$ endowed with the reduced structure is a Lagrangian subvariety. The singular locus $X^s$ does not surject onto $B$.
    \item All smooth fibers are abelian varieties.
    \item If, in addition, $X$ is irreducible symplectic, then $B$ is Fano. Moreover, if $B$ is smooth, then $B \cong \P^n$.
    \end{enumerate}

\section{Preliminary results}

Here we recall a few important results that we are going to use 
in the proofs of our main theorems. 

\hfill

\theorem(Buzzard, Hubbard, \cite{Bu_Hu}) \label{BH2}
Let $V$ be an algebraic subset of codimension two in $\C^N$. 
Then there exists a holomorphic embedding $\C^N \hookrightarrow 
\C^N \setminus V$ avoiding the subset $V$. 

\hfill


%
%
%
%
%
%
%

One of our arguments for our main technical theorem is based on the existence of N\'eron models in codimesion one for Lagrangian fibrations on hyperk\"ahler manifolds.  
%
%
%
%
%
An immediate corollary of this is the existence of an abelian group scheme structure outside a codimension two subset of $B$ over which $\pi$ is a complex analytic torsor in codimension one via the standard facts about the N\'eron model. 

\hfill

\theorem(N\'eron, \cite[Section 1.4]{BLR}) \label{BLR2} Let $S$ be a connected Dedekind domain with field of fractions $K$ and Let $A_K$ be an abelian variety over $K$. Then $A_K$ admits a global N\'eron model over $S$, i.e. a variety $X$ over $S$ that is smooth, separated and of finite type that satisfies the N\'eron mapping property: for each smooth $S$-scheme $Y$ and $K$-morphism $u_K: Y_K\to X_K$, there is a unique $S$-morphism $u:Y\to X$ extending $u_K$.

\hfill

\theorem(Campana, Winkelmann, \cite{CW}) \label{CW2}
Let $X$ be a rationally connected variety.  Then $X$ admits a Zariski dense entire holomorphic curve.

%

\hfill

\theorem \label{logFano}
Log $\mathbb Q$-Fano varieties 
are 
rationally connected. 

\hfill

\proof
\cite{Lu02, QZ}
\endproof

\hfill

We now apply these theorems in the next sections to some of the standard hyperk\"ahler and Calabi-Yau examples to see that they are anti-hyperbolic, in the sense of  admitting dense entire holomorphic curves or dominating holomorphic or meromorphic maps from the complex affine space of the same dimension. In \cite{klv} we established vanishing of the Kobayashi pseudometric for $K3$ surfaces and for many hyperk\"ahler cases admitting Lagrangian fibrations as another measure of anti-hyperbolicity. 

\hfill



\section{Dominability criteria}


\definition
An $n$-dimensional complex manifold $M$ is {\em (meromorphically) 
dominable} by $\C^n$ if there is a meromorphic dominating map 
$F: \C^n \dashrightarrow M$ such that the Jacobian determinant 
$\det(DF)$ is not identically zero in the domain of $F$. If the map $F$ 
is holomorphic, then $M$ is {\em holomorphically dominable} by $\C^n$. 
By abuse of notation, we will at times use dominable to mean dominable by $\C^n$.

\hfill

\example
In \cite{_Bu_Lu_}, it is shown that if $S$ is a compact 
complex surface birationally equivalent to an elliptic K3 surface or to a 
Kummer K3 surface, then $S$ is holomorphically dominable by $\C^2$.

\hfill

\definition
A complex projective variety $M$ is {\em unirational} if there is a dominant 
rational map $\C\P^n \dashrightarrow M$. It is {\em stably rational} if 
$M\times \C\P^l$ is rational for some $l$. It is {\em rationally connected} 
if any two points of M are contained in a rational curve, i.e. a holomorphic 
image of $\C\P^1$. 

\hfill 

\remark 
The following implications are clear (among projective varieties): 
{\em rational $=>$ stably rational $=>$ unirational $=>$ meromorphically 
dominable by $\C^m$; unirational $=>$ rationally connected.} 

\hfill

\proposition
Let $S$ be a compact complex surface dominable by $\C^2$. 
Then the Hilbert scheme $Hilb^n (S)$ is meromorphically dominable by 
$\C^{2n}$. 

\hfill

\proof
Indeed, $\C^{2n}$ dominates $S^n$. 
There is a natural quotient map 
$$S^n \arrow S^n/ \Sigma_n = \Sym^n (S),$$
where $\Sigma_n$ is the symmetric group of $n$ symbols, and therefore 
$\Sym^n(S)$ is holomorphically dominable by $\C^{2n}$. Since $\Hilb^n(S)$ 
is the resolution of singularities of $\Sym^n (S)$, they are birational, 
and hence $\Hilb^n(S)$ is meromorphically dominable by $\C^{2n}$. 
\endproof

\hfill

\proposition
Let $S$ be an elliptic K3 surface, then the Hilbert scheme $Hilb^n(S)$ is 
holomorphically dominable by $\C^{2n}$. 

\hfill

\proof 
Consider an elliptic fibration structure $f: S \arrow B$. Then $B=\C\P^1$ as $q(S)=0$ and the canonical bundle formula applied to $f$ shows that $f$ has no multiple fibres and hence admits local (holomorphic) sections above any point by a classical argument, cf. \cite{LuMM}. Now choose an 
affine open chart $i: \C \subset \C\P^1$. Since $\C$ is Stein,
there is a global holomorphic section $\sigma$ over $\C$ by patching up the local ones via the classical Cousin argument, see \cite{_Bu_Lu_}. 
Hence, there is a holomorphic dominable map from $\C^2$ to $S$ by the following construction: The relative tangent bundle of the fibration restricted to the section $\sigma(\C)$ is a holomorphically trivial bundle $N$ over $\C$ as $\C$ is  Stein and contractible (to a lower dimensional skeleton). Hence $N$ is biholomorphic to $\C^2$. The same argument applied to the restriction of $f$ over $U:=B\setminus \text{disc}(f)$, viewing $f$ there as an abelian torsor over $U$, shows that it is the trivial torsor and hence an abelian scheme over $U$. Hence, the resulting group structure on the smooth fibres of $f$ extends to an open subset of any fiber above $\C$ by the N\'eron mapping property via the N\'eron model, and the relative exponential map from $N$ is defined on all of $N$ and gives the desired dominable map. 

By functoriality of the construction of Hilbert schemes, we have a morphism  
$\Hilb^n(S) \arrow \Hilb^n(\C \P^1) = \C \P^n$.
It has an induced section $\tilde \sigma$ whose image in $\Hilb^n(S)$ is $\Hilb^n (\sigma (\C)) \cong \Hilb^n (\C) = \Sym^n (\C)$. 
Again, the same reasoning and construction via the exponential map defined on the normal bundle to this section shows that 
$\Hilb^n(S)$ is holomorphically dominable. 
\endproof

\hfill

\proposition
Let $A$ be an abelian surface. Then the generalized Kummer variety 
$K^n(A)$ is meromorphically dominable by $\C^{2n}$. 

\hfill

\proof
Indeed, notice that $K^n(A)$ is the kernel of the composition morphism 
$$\Hilb^{n+1}(A) \arrow \Sym^{n+1}(A) \arrow A$$ 
and $A^n$ is the kernel of the 
summation map $\Sigma: A^{n+1} \arrow A$. 
Let $K^n_s(A)$ be the birational image of $K^n(A)$ in $\Sym^{n+1}(A)$. 
Since $K^n_s(A)=A^n/ \Sigma_n$, 
and $A^n$ is dominable by $\C^{2n}$, 
$K^n(A)$ is meromorphically dominable by $\C^{2n}$. 
\endproof

\hfill

\remark
If one assumes the (rational) SYZ conjecture, and if the Picard rank of a hyperk\"ahler 
manifold or of a primitive symplectic variety $M$ is greater than 5, it  would admit a double fibration, induced by two non-proportional parabolic (a.k.a. isotropic) classes in the Picard lattice.\footnote{
Non-propotional means that neither cohomology class is a multiple of the other one.} 
By \cite[Proposition 5.5]{KL}, two such parabolic classes given the said conjecture induce transversal rational Lagrangian fibrations if $b_2(M)\ge 5$; By transversal, we mean that they induce a generically finite map to the product of their base. One can also deduce this without the SYZ ansatz if instead one of the fibrations is almost holomorphic, a situation that one can always reduce to by a birational change of $M$ in the case of projective $M$ with $b_2(M)\ge 5$, and if the two fibrations are not induced by a birational change of the base. 
This observation together 
with \ref{rm_double} that double fibrations are dense in the Teichm\"uller 
space makes it interesting to study dominability of a hyperk\"ahler 
manifold admitting two transversal Lagrangian fibrations. 

\hfill

\theorem \label{main1}
Let $M$ be a hyperk\"ahler manifold or a primitive symplectic variety of dimension $2n$. 
If $M$ admits two rational Lagrangian fibrations $f_1: M \arrow B_1$ 
and $f_2: M \arrow B_2$ that are transversal, i.e. $(f_1,f_2):M\dashrightarrow B_1\times B_2$ is generically finite, then $M$ is meromorphically dominable by $\C^{2n}$. 

\hfill

\proof
Let $F_i$ be a general fiber of $f_i: M \arrow B_i$, for $i=1,2$. 
Then their intersection $F_1 \cap F_2$ is finite, as seen in the proof of 
Theorem 2.11 of \cite{klv}.
So the morphism $\Xi: F_2 \arrow B_1$ is generically finite and in fact finite as $F_2$ is an abelian variety. 
Consider now its fiber product $Z=F_2 \times_{B_1} M$ with $f_1: M \arrow B_1$ and let $\tilde Z \arrow Z$ be its normalization. 
Then the projection $Z\to F_2$ is the pull-back abelian fibration ${\tilde f}$ of $f_1$ via the finite base change $\Xi$. Since $f_1: M \arrow B_1$ is an abelian  
fibration, the projection ${\tilde f} : {\tilde Z} \arrow F_2$ is also 
an abelian fibration but now equipped with a section $\sigma$. 

$$\begin{array}{ccccccc}
{\tilde Z} & \longrightarrow & Z & =
& F_2 \times_{B_1} M  
& \longrightarrow & M \\ 
   &    &      &     &    \downarrow  &               & \downarrow \\
   &    &      &     &    F_2 & \longrightarrow & B_1 \\
\end{array}$$

Let $\pi: \C^n \arrow F_2$ be the universal cover of $F_2$ and let  
$Z' = \C^n \times_{F_2} {\tilde Z}$ be the fiber product of 
$\pi: \C^n \arrow F_2$ and ${\tilde f} : {\tilde Z} \arrow F_2$. 
Denote by $U' \subset Z'$ the locus of all smooth points of 
$f': Z' \arrow \C^n$. 

$$\begin{array}{ccccc}
U' & \hookrightarrow & Z' & \longrightarrow & \tilde Z \\ 
      &     &    \downarrow  &               & \downarrow \\
      &     &    \C^n & \longrightarrow & F_2 \\
\end{array}$$

Since there is a 
section $\sigma$ of ${\tilde f} : {\tilde Z} \arrow F_2$ by construction, this 
gives rise to a  section $\tilde {\sigma}$ of $f': Z' \arrow \C^n$. 
The existence of ${\sigma}$ implies that an open subset of  
$\tilde Z$ proper outside the discriminant locus of $\tilde f$ can be identified outside a codimension two subset $S'$ on the base $F_2$ with the fiberwise dual to the relative Picard variety $\Pic^0=({\tilde f}_{*1} {\cal O}_{\tilde Z})/({\tilde f}_{*1} \Z)$ of $\tilde f$. 
Since 
${\cal F}:=({\tilde f}_{*1} {\cal O}_{\tilde Z})^{\vee}$ is a torsion free sheaf over $F_2$, it is locally free outside a codimension two subset $S\subset F_2$ and it is globally generated after tensoring by an ample invertible sheaf. 
Pulling back this by $\pi$ shows that an open subset of 
$U'$ proper outside the discriminant locus can be identified with $[({f'}_{*1} {\cal O}_{Z'})/({f'}_{*1} \Z)]^{\vee}$ and $\pi^*{\cal F}=({f'}_{*1} {\cal O}_{Z'})^{\vee}$ is globally generated as any invertible sheaf is trivial over $\C^n$ (by appealing to the Oka-Grauert principle for example, $\C^n$ being Stein and homotopic to a point), i.e., there is an exact sequence of sheaves $\oplus_{i=1}^m{\cal O}_{\C^n}\to \pi^*{\cal F}\to 0$. 
Now, in the complement $V$ of the complex codimension two complex analytic subset $\pi^{-1}(S\cup S')$ in $\C^n$, $\pi^*{\cal F}$ is given by a vector bundle that is dominated via this exact sequence by a trivial rank $m$ vector bundle, which can be identified with the analytic subset $V\times \C^m\subset \C^{n+m}$ with complex codimension two complement. By composition, this yields a dominating holomorphic map from $V\times \C^m$ to $U', Z', Z$ and to $M$. Since $M$ is Moishezon, as it has two transveral Lagrangian fibrations (cf. \ref{proj}), it is projective and we may extend the above dominating map across codimension two to yield a meromorphic domination of $M$ by $\C^{n+m}$ and hence also by $\C^{2n}$ on taking a general linear subspace $\C^n\subset \C^{m}$. 
\endproof

\hfill

\remark 
By the results of Kamenova-Lehn \cite[Lemma 2.12]{KL}, it follows that \ref{main1} is valid for rational Lagrangian fibrations that are not induced by a birational change of the base on projective primitive symplectic varieties $M$, assuming $b_2(M) \geq 5$. 
Very minor modifications of  \ref{main1} would be necessary in the singular case, because we are interested in meromorphic dominability. The same apply to \ref{main2}.

\hfill

\remark
We note that, by the same token, unirational fibrations over a unirational 
base are dominable by $\C^m$. This can be seen by restricting to an affine 
chart after a finite base change to a projective space base. In contrast, 
one can construct conic fibrations over $\C\P^2$  that are not stably 
rational, as in \cite{HKT}. 

\hfill

\theorem \label{main2}
Let $M$ be a hyperk\"ahler manifold or a primitive sympletive variety of dimension $2n$. 
If $M$ is projective and admits a Lagrangian fibration $f: M \arrow \C\P^n$ 
with no multiple fibers in codimension one of the base, then $M$ is 
holomorphically dominable by $\C^{2n}$. 

\hfill

\proof
Let $\Delta \subset \C\P^n$ be the discriminant locus of $f$, and let 
$\Delta_m \subset \C\P^n$ be a codimension two subset outside of which there are no multiple fibers. 
Restrict the fibration $f$ to an open affine chart $\C^n \subset \C\P^n$, 
i.e., $f:M_0 \arrow \C^n$. Then $f$ defines an abelian variety over the generic point of $\C^n$. By the Neron mapping property, there is a codimension two algebraic subset $S$ of $\C^n\setminus \Delta_m$ such that a Neron model $\pi: N\to \C^n\setminus \{S\cup\Delta_m\}$ exists for $f$ above $\C^n\setminus \{S\cup\Delta_m\}$.
By \ref{BH2} there is an injective map $i: \C^n \hookrightarrow \C^n \setminus 
\Delta_m$ avoiding $\Delta_m$, i.e., the pull-back abelian fibration 
$f_0: i^*M_0 \arrow \C^n$ doesn't have multiple fibers. Therefore, there are 
local sections of $f_0$ which we can glue together to a global section $\sigma$ by the solution to the Cousin problem as $\C^n$ is Stein. Here by local, we naturally mean the local holomorphic sections in the Euclidian topology. Now, the relative tangent bundle of the fibration restricted to the section $\sigma(\C^n)$ is a holomorphically trivial bundle $N$ over $\C^n$ as $\C^n$ is  Stein and contractible (to a lower dimensional skeleton). Hence, $N$ is biholomorphic to $\C^{2n}$. The same argument applied to $f_0$, viewing it as an abelian torsor over $\C^n$ in the holomorphic category, shows that it is the trivial torsor, and hence an abelian scheme over $\C^n$. Therefore, the resulting group structure on the smooth fibres of $f_0$ extends to an open subset of any fiber above $\C^n$ by the N\'eron mapping property, and the relative exponential map from $N$ is defined on all of $N$ and gives the desired dominable  map.
\endproof

\hfill

\example
O'Grady's exceptional examples $O_6$ and $O_{10}$ are deformation equivalent to Lagrangian fibrations \cite{MO, MR}. They have no multiple fibers in codimension one on the base, and therefore are dominable by $\C^{6}$ and $\C^{10}$, respectively. 
Applying \ref{_dense_ope_Lagra_Theorem_}, we see there are 
densely many Lagrangian fibrations on the deformations of  $O_6$ and $O_{10}$. 
The fibrations without multiple fibers in codimension one are dominable. 

\hfill

\remark
In \cite{klv} the authors together with M. Verbitsky have proved that the 
Kobayashi pseudometric and the infinitesimal Kobayashi pseudometric  of the known examples of hyperk\"ahler manifolds vanish. Here we have shown that some of these examples are dominable by $\C^{2n}$ which is a much stronger property.

\section{Dense entire curves}

Using the main result of Campana-Winkelmann \cite{CW} on the existence of dense entire curves on rationally connected varieties and the fact that log-Fano varieties are rationally connected by \ref{logFano}, we have the following application to anti-hyperbolicity of varieties. 

\hfill

\theorem  \label{main5}
Given a projective manifold $M$ with an abelian fibration $f: M \arrow B$ with no multiple fibers in codimension one over a log-Fano variety $B$, then there exists a dense entire curve.

\hfill

\proof Let $p:\C\to B$ be a holomorphic map with Zariski dense image (cf. \ref{logFano} and \ref{CW2}) and consider its pullback abelian fibration $\tilde f: \tilde M\to \C$. Now by \cite{_Bu_Lu_} we can construct a holomorphic curve that passes through a point of $\tilde M$ an infinite number of times. This implies that its image in $M$ is Zariski dense. 
\endproof

\hfill

{\bf Acknowledgements:} The first named author is grateful to UQAM for their hospitality when starting this project in 2016, and for their hospitality again in 2023. She also thanks the AMS for the MCA 2017 travel award, which made it possible to continue working on the project during the MCA 2017 in Montr\'eal. The authors thank CRM and Universit\'e de Montr\'eal for their hospitality in 2024, where we finalized the main part of this project. We are very grateful to Jason Starr, Xi Chen, Christian Lehn, Claire Voisin and Frederic Campana for their help and interesting discussions. 

\hfill

\noindent {\sc Ljudmila Kamenova\\
Department of Mathematics, 3-115 \\
Stony Brook University \\
Stony Brook, NY 11794-3651, USA,} \\
\tt kamenova@math.sunysb.edu
\\
 
\noindent {\sc Steven Lu\\
D\'epartment de Math\'ematiques, PK-5151\\
Universit\'e du Qu\'ebec \`a Montr\'eal (UQAM)\\
C.P. 8888 Succersale Centreville H3C 3P8,\\ 
Montr\'eal, Qu\'ebec, Canada,} \\
\tt lu.steven@uqam.ca


\begin{thebibliography}{GMP}

\bibitem[Be]{Beauville} Beauville, A., 
        {\em Vari\'et\'es K\"ahleriennes dont la
        premi\`ere classe de Chern est nulle,} 
        J. Diff. Geom.  18 (1983) 755-782. 


\bibitem[BH]{Bu_Hu} Buzzard, G.; Hubbard, J., 
        {\em A Fatou-Bieberbach domain avoiding a neighborhood of a 
         variety of codimension 2,} Math. Ann. 316 (2000) 699-702. 

\bibitem[BL1]{_Bu_Lu_} Buzzard, G., Lu, S., 
        {\em Algebraic Surfaces Holomorphically Dominable by $\C^2$}, 
        Invent. Math.  139 (2000) 617-659. 

\bibitem[BL2]{BL1} Buzzard, G., Lu, S., 
        {\em Double sections, dominating maps, and the Jacobian fibration}, 
       American Journal of Mathematics 122, No. 5 (2000) 1061-1084.
       
\bibitem[BLR]{BLR}
Bosch, S.,  L\"utkebohmert, W.,  Raynaud M.,
{\em N\'eron Models},
Ergeb. Math. Grenzgeb. 21, Springer, Berlin, 1990.


\bibitem[Cam]{Camp04} Campana, F.,
 {\em Special varieties and classification theory,}
  Ann. Inst. Fourier (Grenoble), 54 (2004) 499-630. 


\bibitem[CW]{CW} Campana, F., Winkelmann, J., 
{\em Dense entire curves in rationally connected manifolds (with an appendix by 
J\'anos Koll\'ar),} 
Algebraic Geometry 10 (2023) 521-553. 

\bibitem[Cat]{_Catanese:moduli_}
Catanese, F., 
{\em A Superficial Working Guide to Deformations and Moduli}, 
arXiv:1106.1368, 
Advanced Lectures in Mathematics, Volume XXVI 
Handbook of Moduli, Volume III, page 161-216 
(International Press).

\bibitem[F]{_Fujiki:HK_}  
Fujiki, A., {\em On the de Rham Cohomology Group of a Compact 
K\"ahler Symplectic Manifold}, Adv. Stud.
Pure Math. 10 (1987) 105-165.

\bibitem[G]{Grassi}
Grassi, A., {\em On minimal models of elliptic threefolds,} 
Math. Ann. 290 (1991) 287-301. 

\bibitem[GL]{GL} Greb, D., Lehn, C., 
{\em Base manifolds for Lagrangian fibrations on hyperk\"ahler manifolds,} 
Int. Math. Res. Not. 19 (2014) 5483-5487. 

\bibitem[HKT]{HKT} Hassett, B., Kresch, A., Tschinkel, Y., 
{\em Stable rationality and conic bundles,} 
Math. Ann. 365 (2016) 1201-1217. 

\bibitem[Hu1]{_Huybrechts:basic_}
Huybrechts, D., 
{\em Compact hyper-K\"ahler manifolds: basic results}, 
Invent. Math. 135 (1999), no. 1, 63-113. Erratum: 
Invent. Math. 152 (2003), no. 1, 209-212. 

\bibitem[Hu2]{H_K3} Huybrechts, D., 
{\em Lectures on K3 surfaces}, Cambridge University Press (2016). 

\bibitem[HX]{HX} Huybrechts, D., Xu, C., 
{\em Lagrangian fibrations of hyperk\"ahler fourfolds,} 
J. Inst. Math. Jussieu 21 (2022) 921-932. 

\bibitem[Hw]{Hwang} Hwang, J.-M., 
{\em Base manifolds for fibrations of projective irreducible symplectic manifolds,} Invent. Math. 174 (2008) 625-644. 


\bibitem[KL]{KL} Kamenova, L., Lehn, C., 
{\em Non-hyperbolicity of holomorphic symplectic varieties,} arXiv:2212.11411 [math.AG].

\bibitem[KV]{_Kamenova_Verbitsky:fibrations_}
Kamenova, L., Verbitsky, M., 
{\em Families of Lagrangian fibrations on hyperkaehler manifolds}, 
Adv. Math. 260 (2014) 401-413. 

\bibitem[KLV]{klv}
Kamenova, L., Lu, S., Verbitsky, M., 
{\em Kobayashi pseudometric on hyperkahler manifolds,} 
J. London Math. Soc. (2014) 90 (2): 436-450. 

\bibitem[Kol]{Kol}
Koll\'ar, J., 
{\em Kodaira’s canonical bundle formula and adjunction,} In Flips for 3-folds and 4-folds, volume 35 of Oxford Lecture Ser. Math. Appl., pages 134–162. Oxford Univ. Press, Oxford, 2007.

\bibitem[KM]{KM}
Koll\'ar, J., Mori, S., 
{\em Birational Geometry of Algebraic Varieties,} 
Cambridge University Press, 1998.

\bibitem[L1]{LuMM} Lu, S., 
{\em Multiply marked Riemann surfaces and the Kobayashi pseudometric on algebraic manifolds,}
Preprint 1999.


\bibitem[L2]{Lu02} Lu, S., 
{\em A refined Kodaira dimension and its canonical fibration,} 
Preprint: arXiv:math.AG/0211029, 2002.




\bibitem[M]{Matsushita} Matsushita, D., 
{\em On fibre space structures of a projective irreducible symplectic manifold,} Topology 38 (1999) 79-83; addendum ibid. 40 (2001) 431-432. 


\bibitem[MO]{MO} Mongardi, G., 
Onorati, C., 
{\em Birational geometry of irreducible holomorphic symplectic tenfolds of O’Grady type,}
Math. Zeit. 300 (2022) 3497 - 3526.

\bibitem[MR]{MR} Mongardi, G., Rapagnetta, A., 
{\em Monodromy and birational geometry of O’Grady’s sixfolds,} J. Math. Pures Appl. 146 (2021) 31 - 68. 

\bibitem[Nak]{Nak}
Nakayama, N., 
{\em The singularity of the canonical model of compact K\"ahler manifolds,} Math. Ann., 280(3) (1988) 509-512.

\bibitem[Nam]{Namba} Namba, M., 
{\em Branched coverings and algebraic functions,} Pitman Research Notes in Mathematics Series, 161 (1987). 

\bibitem[O1]{og1} O'Grady, K., 
        {\em Desingularized moduli spaces of sheaves on a K3,} 
        J. Reine Angew. Math. 512 (1999) 49 - 117. 

\bibitem[O2]{og2} O'Grady, K., 
        {\em A new six-dimensional irreducible symplectic variety,} 
        J. Algebr. Geom. 12 (2003) 435 - 505. 

\bibitem[Og]{Oguiso} Oguiso, K., {\em On algebraic fiber space structures on Calabi-Yau threefolds,} International Journal of Mathematics, Vol. 4 (1993) 439-465. 

\bibitem[Ou]{Ou} Ou, W., 
{\em Lagrangian fibrations on symplectic fourfolds,} 
J. Reine Angew. Math. 746 (2019) 117-147. 


\bibitem[R]{Reid} Reid, M., 
{\em Hyperelliptic linear systems on a K3 surface,} J. London Math. Soc. (2) 13 (1976) 427-437.

\bibitem[Ulu]{Ulu}
Uluda\u{g}, A.M., {\em Orbifolds and Their Uniformization,} In: Holzapfel, RP., Uluda\u{g}, A.M., Yoshida, M. (eds) Arithmetic and Geometry Around Hypergeometric Functions. Progress in Mathematics, vol 260. Birkh\"auser Basel (2007).



\bibitem[V]{_Verbitsky:ergodic_}
Verbitsky, M., 
{\em Ergodic complex structures on hyperk\"ahler manifolds,} 
Acta Math. 215 (2015) 161-182. 

\bibitem[Z]{QZ} Zhang, Q., 
{\em Rational connectedness of log 
Q-Fano varieties,} 
J. Reine Angew. Math. 590 (2006) 131 - 142.

\end{thebibliography}
\end{document}